\documentclass{article}

\usepackage{amsmath,amssymb,amsfonts,amsthm}

\newtheorem{theorem}{Theorem}[section]
\newtheorem{proposition}[theorem]{Proposition}
\newtheorem{lemma}[theorem]{Lemma}

\numberwithin{equation}{section}

\theoremstyle{remark}
\newtheorem{remark}[theorem]{Remark}

\newcommand{\II}{\mathop{\mathrm{II}}\nolimits}
\newcommand{\End}{\mathop{\mathrm{End}}}
\newcommand{\Ric}{\mathop{\mathrm{Ric}}\nolimits}
\newcommand{\Mat}{\mathop{\mathrm{Mat}}\nolimits}
\newcommand{\Id}{\mathop{\mathrm{Id}}\nolimits}
\renewcommand{\Pr}{\mathop{\mathrm{Pr}}\nolimits}

\author{Artem Pulemotov \\ \mbox{} \\ Department of Mathematics \\ The University of Chicago \\
5734~S.~University~Ave., Chicago,~IL 60637-1514, USA \\
Email: \texttt{artem@math.uchicago.edu}}

\title{Quasilinear parabolic equations and the Ricci flow on
manifolds with boundary}

\date{}

\begin{document}

\maketitle

\begin{center}

\end{center}

\begin{abstract}
The first part of the paper discusses a second-order quasilinear
parabolic equation in a vector bundle over a compact manifold $M$
with boundary $\partial M$. We establish a short-time existence
theorem for this equation. The second part of the paper is devoted
to the investigation of the Ricci flow on $M$. We propose a new
boundary condition for the flow and prove two short-time existence
results.

\end{abstract}

\section{Introduction}

The present paper is motivated by the desire to investigate
geometric evolutions on a compact manifold~$M$ with boundary
$\partial M$. Our first goal is to study a second-order quasilinear
parabolic equation in a vector bundle over $M$. We then apply the
obtained results to the analysis of the Ricci flow on $M$. Let us
explain the essence and the history of the problems to be
considered.

A significant step in the investigation of geometric evolutions on
$M$ is to acquire information about second-order quasilinear
parabolic equations in vector bundles over $M$. Particularly, it is
important to have a short-time existence theorem that would cover a
wide range of boundary conditions and produce a solution with ample
differentiability properties. Until now, such a theorem did not
appear in the literature. Even in the case where~$M$ is the closure
of a domain in $\mathbb R^n$, there was no published result that
would meet the demands of the applications to geometric evolutions.
In Section~\ref{sec_parab} of the present paper, we make an effort
to remedy this situation. We establish a short-time existence
theorem for a second-order quasilinear parabolic equation in a
vector bundle over $M$. No assumptions are imposed on the geometry
of $M$. Yet even in the the case where~$M$ is the closure of a
domain in $\mathbb R^n$, our result is somewhat different from the
results that previously appeared in print. We will now describe it
in more detail.

Fix a Riemannian metric on $M$. Let $E$ be a vector bundle over $M$.
Suppose $E$ is equipped with a fiber metric and a connection
$\nabla$. We focus on the equation
\begin{align}\label{intro_eq}
\frac\partial{\partial t}u(x,t)-H^{ij}(u(x,t),t)\nabla_i\nabla_j
u(x,t)=F(u(x,t),\nabla u(x,t),t)
\end{align}
for a section $u$ of $E$ depending on $t\ge0$. Here, $H$ is a smooth
map from $E\times[0,\infty)$ to the space of symmetric (2,0)-tensors
over $M$, and $F$ is a smooth map from $E\times(T^*M\otimes
E)\times[0,\infty)$ to $E$. The meaning of the rest of the notation
should be easy to infer from the context. In the beginning of
Section~\ref{sec_parab}, we explain it pedantically. Suppose now
that $E_{\partial M}$ is the restriction of the bundle $E$ to
$\partial M$ and $W$ is a subbundle of~$E_{\partial M}$. We
supplement equation~\eqref{intro_eq} with the boundary conditions
\begin{align}\label{intro_BC}
\Pr_Wu(x,t)&=o(x), \notag \\
\Pr_{W^\bot}\left(H^{ij}(u(x,t),t)\nu_i(x)\nabla_ju(x,t)\right)&=\Psi(u(x,t),t).
\end{align}
Here, $o$ is the zero section of $E$, and $\nu$ is the outward unit
normal covector field on $\partial M$. The smooth map $\Psi$ acts
from $E_{\partial M}\times[0,\infty)$ to $W^\bot$. The first line
in~\eqref{intro_BC} should be thought of as the Dirichlet boundary
condition. It is imposed on $u$ inside $W$. The second line
in~\eqref{intro_BC} may be looked at as a nonlinear nonhomogeneous
Neumann condition. It is imposed inside $W^\bot$. Finally, we
supplement~\eqref{intro_eq} with the initial condition
\begin{align}\label{intro_IC}
u(x,0)=u_0(x).
\end{align}
In this formula, $u_0$ is a smooth section of $E$. The main result
of Section~\ref{sec_parab} requires two additional assumptions.
First, we demand that equation~\eqref{intro_eq} be parabolic.
Second, we impose a compatibility condition near $\partial M$ when
$t=0$. If these assumptions are satisfied, the main result of
Section~\ref{sec_parab} tells us that
problem~\eqref{intro_eq}--\eqref{intro_BC}--\eqref{intro_IC} has a
solution on $M\times[0,T)$ for some $T>0$. This result appears below
as Theorem~\ref{thm_parab_ex}. We point out that the solution it
produces is smooth on $M\times(0,T)$. A problem akin
to~\eqref{intro_eq}--\eqref{intro_BC}--\eqref{intro_IC} was studied
in W.-X.~Shi's paper~\cite{WXS89}. That work, however, only allowed
Dirichlet-type boundary conditions. It turns out that the
nonlinearities in the second line of~\eqref{intro_BC} contribute
substantially to the difficulty of the question of the existence of
solutions to~\eqref{intro_eq}--\eqref{intro_BC}--\eqref{intro_IC}.

Let us assume for a moment that the manifold $M$ is the closure of a
domain in $\mathbb R^n$ and $E$ is the product $M\times\mathbb R^d$
carrying the standard fiber metric and connection. Note that, even
under these assumptions, it is possible for $W$ to be a nontrivial
bundle. We encounter such a phenomenon when dealing with the Ricci
flow later in the present paper. Suppose for now, however, that
$W=\partial M\times\mathcal R$ with $\mathcal R$ being the space of
all $(e_1,\ldots,e_d)\in\mathbb R^d$ such that
$e_{d'+1}=\ldots=e_d=0$ for a fixed $d'$ between 0 and $d$. In this
case, equation~\eqref{intro_eq} may be viewed as a second-order
quasilinear system. The section $u$ takes the form
$u=(u^1,\ldots,u^d)$ with each $u^m$ real-valued.
Formulas~\eqref{intro_BC} become Dirichlet boundary conditions on
$u^1,\ldots,u^{d'}$ and nonlinear nonhomogeneous Neumann conditions
on $u^{d'+1},\ldots,u^d$. In the case where $M$ is the closure of a
domain in $\mathbb R^n$, $E=M\times\mathbb R^d$, and $W=\partial
M\times\mathcal R$, problems similar
to~\eqref{intro_eq}--\eqref{intro_BC}--\eqref{intro_IC} were
extensively studied. Several different methods were proposed to
prove the existence of solutions. For example, the
papers~\cite{HA86,HA90} by H.~Amann used abstract
functional-analytic techniques. The work~\cite{PW91} by
P.~Weidemaier employed a more straightforward fixed-point argument
in a Sobolev-type space. The reader should
see~\cite{OLVSNU68,MGGM87,PABT88} for other approaches.

No major restrictions are imposed in Section~\ref{sec_parab} on the
geometry of $M$, $E$, or $W$. But even in the case where $M$ is the
closure of a domain in $\mathbb R^n$, $E=M\times\mathbb R^d$, and
$W=\partial M\times\mathcal R$, the material we present there did
not previously appear in the literature in the same form. When we
restrict our attention to this geometrically trivial situation, the
theorem in the introduction of H.~Amann's paper~\cite{HA90} is
somewhat similar to our Theorem~\ref{thm_parab_ex}. However, that
result involves a different compatibility condition. We should
remark that the arguments in~\cite{HA90} are rather complicated. It
seems that adapting them to the setting of manifolds and vector
bundles would be a tedious task. When $M$ is the closure of a domain
in $\mathbb R^n$, $E=M\times\mathbb R^d$, and $W=\partial
M\times\mathcal R$, the reasoning in P.~Weidemaier's
paper~\cite{PW91} is akin to much of our reasoning in
Section~\ref{sec_parab}. On the other hand, the work~\cite{PW91} is
concerned with a narrower range of boundary conditions. Besides, it
does not touch upon the issue of the smoothness of solutions.

As we previously declared, our investigation of second-order
quasilinear parabolic equations was motivated by the desire to study
geometric evolutions on manifolds with boundary. Let us say a few
words about the specific applications of Theorem~\ref{thm_parab_ex}.
In Section~\ref{sec_RF}, we establish two short-time existence
results for the Ricci flow on a manifold with boundary.
Theorem~\ref{thm_parab_ex} is a crucial ingredient in our
considerations. The paper~\cite{NCLG10} studies the Yang-Mills heat
flow on 2- and 3-dimensional manifolds with boundary and utilizes it
for the purposes of quantum field theory. Among other things, the
authors of~\cite{NCLG10} prove the existence of solutions to the
flow; see also~\cite{AP08}. We speculate that
Theorem~\ref{thm_parab_ex} can be used to simplify their arguments.
More applications of this nature may emerge in the near future.
Theorem~\ref{thm_parab_ex} seems to be a convenient tool for proving
the existence of solutions to various geometric evolutions.

Section~\ref{sec_RF} of the present paper focuses on the Ricci flow
on the manifold $M$. Our goal is to introduce a new boundary
condition for the flow and establish two short-time existence
results. More precisely, consider the equation
\begin{align}\label{intro_RF}
\frac{\partial}{\partial t}g(x,t)=-2\Ric^g(x,t)
\end{align}
for a Riemannian metric $g$ on $M$ depending on the parameter
$t\ge0$. We supplement~\eqref{intro_RF} with the initial condition
\begin{align}\label{Intro_RiIC}
g(x,0)=\hat g(x).
\end{align}
Here, $\hat g$ is a smooth Riemannian metric on $M$.
Equation~\eqref{intro_RF} is the Ricci flow equation on $M$. To
learn about its history, intuitive meaning, technical peculiarities,
and geometric applications, the reader should refer to the many
quality books on the subject, such
as~\cite{BCDK04,BCPLLN06,PT06,JMGT07}. Examples of how it comes up
in mathematical physics may be found
in~\cite{MHTW06,TOVSEW06,GHTSCW07,TOEW07} and other papers. One more
interesting application is to the regularization of non-smooth
Riemannian metrics; see, e.g.,~\cite{MS02,MS05}.

The Ricci flow on manifolds with boundary is not yet deeply
understood. The root of all evil lies in the fact that
equation~\eqref{intro_RF} is not parabolic. For this reason, it is
difficult to find geometrically meaningful and analytically
appealing boundary conditions for solutions
of~\eqref{intro_RF}--\eqref{Intro_RiIC}. It would be natural to
demand, for example, that the metric induced by $g$ on $\partial M$
always coincide with the metric induced by $\hat g$. But so far,
nobody knows how to prove the short-time existence of solutions
under such a requirement. Progress towards finding boundary
conditions to go with~\eqref{intro_RF}--\eqref{Intro_RiIC} was made
by Y.~Shen in his dissertation~\cite{YS92}. Those results were also
published in the paper~\cite{YS96}. Y.~Shen considered the case
where the second fundamental form $\hat\II$ of the boundary with
respect to $\hat g$ satisfied the equality $\hat\II(x)=\tau\hat
g(x)$ with $\tau\in\mathbb R$ for all $x\in\partial M$. In other
words, he assumed that $\partial M$ was umbilic\footnote{There is
ambiguity in the literature as to the use of the term ``umbilic" in
this context. See the discussion in~\cite{JC09}.} when $t=0$. He was
then able to prove the existence of $T>0$ and a solution $g$ to
problem~\eqref{intro_RF}--\eqref{Intro_RiIC} on $M\times[0,T)$ such
that the second fundamental form $\II$ of $\partial M$ with respect
to $g$ satisfied the equality
\begin{align}\label{intro_umbilic}
\II(x,t)=\tau g(x,t)
\end{align}
for all $x\in\partial M$ and $t\in[0,T)$. Not much is known about
the behavior of $g$ for large $t$. This question was addressed
in~\cite{JC09} under additional assumptions. We should point out
that the case where $\tau=0$ is somewhat special. If $\tau=0$, then
$\partial M$ is totally geodesic with respect to $\hat g$. In this
situation, one can say a few things about how the solution $g$
produced in~\cite{YS92} behaves for large $t$;
see~\cite{YS92,YS96,XCTD06}.

No new boundary conditions for the Ricci flow in dimensions higher
than 2 have been proposed in the literature since the publication of
Y.~Shen's dissertation. A certain amount of work, however, has been
done on surfaces. The list of relevant texts
includes~\cite{TL93,SB02,JC07}. While not much is known today about
the Ricci flow on manifolds with boundary, it is clear that results
in this area would have significant geometric applications. They
would also be useful to mathematical physicists;
see~\cite{MHTW06,GHTSCW07}.

In section~\ref{sec_RF}, we propose a new boundary condition on the
solutions of the Ricci flow and prove two short-time existence
results, a theorem and a proposition. Let $\hat{\mathcal H}$ be the
mean curvature of $\partial M$ with respect to $\hat g$. Our theorem
assumes that $\hat{\mathcal H}$ is equal to a constant $\mathcal
H_0\in\mathbb R$ everywhere on $\partial M$. It then claims that
there exist $T>0$ and a solution $g$ to
problem~\eqref{intro_RF}--\eqref{Intro_RiIC} on $M\times[0,T)$ such
that the boundary condition
\begin{align*}
\mathcal H(x,t)=\mu(t)\mathcal H_0
\end{align*}
holds for all $x\in\partial M$ and $t\in[0,T)$. Here, $\mu$ is a
function that may be thought of as a normalization factor, and
$\mathcal H$ is the mean curvature of $\partial M$ with respect to
$g$. Our proposition touches upon the question of the behavior of
the Ricci flow on manifolds with convex boundary. This question is
natural, and it is related to some of the material
in~\cite{YS92,YS96,JC09,MBXCAP10}. Among other things, our
proposition implies that, if $\partial M$ is convex with respect to
$\hat g$, we can find $T>0$ and a solution $g$
of~\eqref{intro_RF}--\eqref{Intro_RiIC} on $M\times[0,T)$ such that
$\partial M$ remains convex with respect to $g$.

The results in Section~\ref{sec_RF} constitute a step towards
understanding the Ricci flow on manifolds with boundary. We suspect
they can also be used for the purposes of regularizing non-smooth
Riemannian metrics on such manifolds. It is worth mentioning that
the proofs of the results in Section~\ref{sec_RF} are based on the
method commonly known as DeTurck's trick. These proofs rely
substantially on Theorem~\ref{thm_parab_ex}.

\section{Parabolic equations in vector bundles}
\label{sec_parab}

Suppose $M$ is a smooth $n$-dimensional ($n\ge2$) manifold with
boundary. We assume that $M$ is compact, connected, and oriented.
The notations $M^\circ$ and $\partial M$ will be used for the
interior and the boundary of~$M$. Consider a smooth vector bundle
$E$ over $M$ with projection $\pi$ and standard fiber $\mathbb R^d$.
We will discuss second-order quasilinear parabolic equations for
sections of $E$ subject to nonlinear nonhomogeneous boundary
conditions. Our goal will be to establish the short-time existence
of solutions. This result will later help us investigate the Ricci
flow.

Throughout Section~\ref{sec_parab}, we fix a smooth Riemannian
metric on $M$. The tangent bundle $TM$ is equipped with the
Levi-Civita connection. The letter $\nu$ will stand for the outward
unit normal covector field on $\partial M$. Let us also fix a smooth
fiber metric in $E$ and a smooth connection in $E$ compatible with
this metric.

Our arguments will involve tensor products of the form $\mathfrak
E=\mathfrak E_1\otimes\cdots\otimes\mathfrak E_k$ with $\mathfrak
E_i$ equal to $TM$, $T^*M$, or $E$ for each $i=1,\ldots,k$. Here,
$T^*M$ designates the cotangent bundle, and $k$ is a natural number.
The Riemannian metric on $M$ and the fiber metric in $E$ generate a
fiber metric in every such $\mathfrak E$. Given $\eta\in\mathfrak
E$, we write $|\eta|$ for its norm. The Levi-Civita connection in
$TM$ and the fixed connection in $E$ give rise to a connection in
$\mathfrak E$. We write $\nabla f$ for the covariant derivative of a
section $f$ of $\mathfrak E$. Our considerations will also involve
second-order differential operators. In particular, $\nabla\nabla f$
stands for the second covariant derivative of~$f$.

We will sometimes employ local coordinates on the manifold $M$. Let
us introduce the corresponding notation. In what follows, we
implicitly assume that a coordinate system $\{x_1,\ldots,x_n\}$ is
chosen in a neighborhood of every point $x\in M$. If $T$ is a
$(k,l)$-tensor at $x$, we write $T_{j_1\ldots j_l}^{i_1\ldots i_k}$
for its components in this coordinate system. Given a section $f$ of
the bundle $\mathfrak E$, the notation $\nabla_if$ stands for its
covariant derivative in the direction of $\frac\partial{\partial
x_i}$. Analogous shorthand is used for second-order operators.
Namely, $\nabla_i\nabla_jf$ means $\nabla\nabla f$ applied to
$\frac\partial{\partial x_i}$ and $\frac\partial{\partial x_j}$. The
Einstein summation convention is in effect.

\subsection{Formulation of the existence theorem}
\label{subsec_parab_form}

Our purpose is to study the solvability of second-order quasilinear
parabolic equations for sections of $E$ subject to nonlinear
nonhomogeneous boundary conditions. To begin with, consider a smooth
mapping
\begin{align*}
H:E\times[0,\infty)\to TM\otimes TM.
\end{align*}
We assume that $H(\eta,t)$ is a symmetric tensor over $\pi(\eta)$
for all $\eta$ and $t$. It will be necessary to impose one more
requirement on~$H$, but we postpone this until a little later.
Meanwhile, consider another smooth mapping
\begin{align*}
F:E\times(T^*M\otimes E)\times[0,\infty)\to E.
\end{align*}
We demand that $\pi(F(\eta,\theta,t))=\pi(\eta)$ for all values of
$\eta$, $\theta$, and $t$. Our attention will be focused on the
equation
\begin{align}\label{parab_eq}
\frac\partial{\partial t}u(x,t)-H^{ij}(u(x,t),t)\nabla_i\nabla_j
u(x,t)=F(u(x,t),\nabla u(x,t),t),\qquad x\in M^\circ,~t\in(0,T),
\end{align}
for a section $u$ of $E$ depending on the parameter $t\in[0,T)$ with
$T>0$. The first step is to supplement~\eqref{parab_eq} with
boundary conditions.

Let $E_{\partial M}$ denote the set of all $\eta\in E$ such that
$\pi(\eta)\in\partial M$. This set has the structure of a vector
bundle over $\partial M$ induced by the structure of $E$. Also,
$E_{\partial M}$ inherits the fiber metric from $E$. Suppose $W$ is
a subbundle of $E_{\partial M}$. Let $W^\bot$ be the orthogonal
complement of $W$ in $E_{\partial M}$. Introduce a smooth mapping
\begin{align*}
\Psi:E_{\partial M}\times[0,\infty)\to W^\bot.
\end{align*}
It is assumed that $\pi(\Psi(\eta,t))=\pi(\eta)$ for all values of
$\eta$ and $t$. We impose the boundary conditions
\begin{align}\label{parab_BC}
\Pr_Wu(x,t)&=o(x), \notag \\
\Pr_{W^\bot}\left(H^{ij}(u(x,t),t)\nu_i(x)\nabla_ju(x,t)\right)&=\Psi(u(x,t),t),\qquad
x\in\partial M,~t\in(0,T),
\end{align}
on the solutions of~\eqref{parab_eq}. Here, $\Pr_W$ and
$\Pr_{W^\bot}$ stand for the orthogonal projections in $E_{\partial
M}$ onto $W$ and~$W^\bot$. The letter $o$ refers to the zero section
of $E$. In essence, we impose the Dirichlet boundary condition on
$u$ inside~$W$ and a nonlinear nonhomogeneous Neumann condition
inside $W^\bot$. To understand what formulas~\eqref{parab_BC} mean
in a geometrically trivial case, the reader may revisit the
introduction to the present paper.

Suppose $u_0$ is a smooth section of $E$. We supplement
equation~\eqref{parab_eq} with the initial condition
\begin{align}\label{parab_InitC}
u(x,0)=u_0(x),\qquad x\in M.
\end{align}
Our goal is to establish the solvability of
problem~\eqref{parab_eq}--\eqref{parab_BC}--\eqref{parab_InitC}. In
order to do so, we need two additional assumptions. The first one is
a parabolicity condition on equation~\eqref{parab_eq}. We suppose
there is a constant $c_1>0$ such that the inequality
\begin{align}\label{parabolicity_cond}
H^{ij}(\eta,t)\xi_i\xi_j\ge c_1|\xi|^2
\end{align}
holds for every $\eta\in E$, $t\in[0,\infty)$, and $\xi\in T^*M$
projecting on $\pi(\eta)$. The reader should see,
e.g.,~\cite[Chapter~4]{PT06} for an elaborate discussion of the
concept of parabolicity in the framework of vector bundles. The
second assumption is the natural compatibility condition
\begin{align}\label{parab_compat}
\Pr_Wu_0(x)&=o(x), \notag \\
\Pr_{W^\bot}\left(H^{ij}(u_0(x),0)\nu_i(x)\nabla_ju_0(x)\right)&=\Psi(u_0(x),0),\qquad
x\in\partial M.
\end{align}
It is now time to state the main result of Section~\ref{sec_parab}.
This result is an existence theorem for
problem~\eqref{parab_eq}--\eqref{parab_BC}--\eqref{parab_InitC}. It
will be utilized in Section~\ref{sec_RF} when we investigate the
Ricci flow.

\begin{theorem}\label{thm_parab_ex}
Consider the initial-boundary value
problem~\eqref{parab_eq}--\eqref{parab_BC}--\eqref{parab_InitC}.
Suppose the parabolicity condition~\eqref{parabolicity_cond} and the
compatibility condition~\eqref{parab_compat} are satisfied. Then
there exist a number $T>0$ and a map $u:M\times[0,T)\to E$ such that
the following requirements are met:
\begin{enumerate}
\item
The equality $\pi(u(x,t))=x$ holds for every $x\in M$ and
$t\in[0,T)$. In other words, $u$ is a section of $E$ depending on
$t\in[0,T)$.
\item
The map $u$ and the covariant derivative $\nabla u$ are continuous
on $M\times[0,T)$. Furthermore, $u$ is smooth on $M\times(0,T)$.
\item
Equalities~\eqref{parab_eq}, \eqref{parab_BC},
and~\eqref{parab_InitC} hold for $u$.
\end{enumerate}
\end{theorem}

\begin{remark}
The number $T>0$ whose existence the theorem asserts is dependent on
the mappings $H$, $F$, $\Psi$, and $u_0$. It may also be affected,
for example, by the Riemannian metric on $M$.
\end{remark}

\begin{remark}\label{rem_restr_HFg}
Suppose $\Xi$ is an open neighborhood of the set $\{u_0(x)\,|\,x\in
M\}$ in $E$. It is not difficult to verify that the theorem still
holds if the mappings $H$, $F$, and $\Psi$ are only defined on
$\Xi\times[0,\infty)$, $\Xi\times(T^*M\otimes E)\times[0,\infty)$,
and $(\Xi\cap E_{\partial M})\times[0,\infty)$. Furthermore, let
$\Omega$ be the set of $(\eta,\theta,t)\in\Xi\times(T^*M\otimes
E)\times[0,\infty)$ such that $\eta$ and $\theta$ project onto the
same point in $M$. The theorem prevails if $F$ is only defined on
$\Omega$.
\end{remark}
\begin{remark}\label{rem_par_comp}
It may be possible to improve the regularity of $u$ on
$M\times[0,T)$ by imposing higher-order compatibility conditions
along with~\eqref{parab_compat};
cf.~\cite[pages~319--321]{OLVSNU68}. But this issue remains beyond
the scope of the present paper.
\end{remark}

\begin{remark}
One may be able to prove an analogue of Theorem~\ref{thm_parab_ex}
in a more general setting. Namely, suppose $\End E$ is the bundle of
endomorphisms of $E$ and the mapping $H$ acts from
$E\times[0,\infty)$ to $TM\otimes TM\otimes\End E$ instead of
$TM\otimes TM$. It is clear how
problem~\eqref{parab_eq}--\eqref{parab_BC}--\eqref{parab_InitC}
should be modified in this case. It may then be possible to
establish an existence result analogous to
Theorem~\ref{thm_parab_ex}. But we do not concern ourselves with
this in the present paper. Let us just mention that the
references~\cite{VS65,HA86,MGGM87,BT89,HA90,PW91} might be helpful.
\end{remark}

Before we can prove Theorem~\ref{thm_parab_ex}, we need to introduce
additional notation, make a few comments, and state a lemma. This
will be done in Sections~\ref{subsec_Holder_sp}
and~\ref{subsec_lin_par}. When the preparations are finished, we
will use a fixed-point argument to produce $T$ and $u$. The last
step will be to establish the smoothness of $u$ by localizing our
equation and appealing to some classical facts from~\cite{OLVSNU68}.

\subsection{Spaces of vector bundle sections}
\label{subsec_Holder_sp}

We will deal with a multitude of spaces of vector bundle sections.
Although some of these spaces are rather classical, they can be
approached from several different viewpoints. In order to exclude
ambiguity, and for the convenience of the reader, we will outline
the definitions with which we will work in this paper.

Let us use the notation $\mathbb R^n_{+,0}$ for the open half-space
$\{(y_1,\ldots,y_n)\in\mathbb R^n\,|\,y_n>0\}$ and the notation
$\mathbb R^n_+$ for the closed half-space
$\{(y_1,\ldots,y_n)\in\mathbb R^n\,|\,y_n\ge0\}$. Fix a real number
$I\in(0,1)$ and an integer number $q>n+2$. We will encounter the
classical Sobolev-type spaces $W^{2,1}_q(\mathbb R^n\times(0,I))$
and $W^{2,1}_q(\mathbb R_{+,0}^n\times(0,I))$ of real-values
functions. Their precise definitions can be found in several sources
such as, for example,~\cite[Chapter~I]{OLVSNU68}. Given a domain
$\Theta$ in $\mathbb R^n$ or $\mathbb R^n_+$ and a number
$\lambda>0$, we will deal with the H\"older-type space
$H^{\lambda,\frac\lambda2}(\Theta\times(0,I))$ of real-valued
functions. Again, one may find its definition
in~\cite[Chapter~I]{OLVSNU68}.

Like in the beginning of Section~\ref{sec_parab}, consider a vector
bundle $\mathfrak E=\mathfrak E_1\otimes\cdots\otimes\mathfrak E_k$.
Let $C^\infty(M\times[0,I]\to\mathfrak E)$ be the set of all the
smooth mappings $\phi:M\times[0,I]\to\mathfrak E$ such that the
projection of $\phi(x,t)$ onto $M$ always equals $x$. Suppose $dx$
and $dt$ are the Riemannian volume measure on $M$ and the Lebesgue
measure on $[0,I]$. We will encounter the space
$L^q(M\times[0,I]\to\mathfrak E,dx\,dt)$. It is the completion of
$C^\infty(M\times[0,I]\to\mathfrak E)$ in the norm
\begin{align*}
\|\phi\|_{L^q(M\times[0,I]\to\mathfrak
E,dx\,dt)}=\bigg(\int_{M\times[0,I]}|\phi|^q\,dx\,dt\bigg)^{\frac1q}.
\end{align*}
Our further arguments will involve the spaces
$L^q(M\times[0,I]\to\mathfrak E,dx\,dt)$ for several different
bundles $\mathfrak E$. It will be convenient to use the same short
notation $L^q_I$ for all these spaces. The norms
$\|\cdot\|_{L^q(M\times[0,I]\to\mathfrak E,dx\,dt)}$ will all be
written as $\|\cdot\|_{L^q_I}$.

Suppose $C^\infty_W(M\times[0,I]\to E)$ is the set of all the smooth
mappings $\phi\in C^\infty(M\times[0,I]\to E)$ such that the
equalities
\begin{align*}
\Pr_W(\phi(x,t))&=o(x),\qquad x\in\partial M,~t\in[0,I],\\
\phi(x,0)&=o(x),\qquad x\in M,
\end{align*}
hold true. Let $\mathcal W_q^I$ stand for the completion of
$C^\infty_W(M\times[0,I]\to E)$ in the norm
\begin{align*}
\|\phi\|_{\mathcal
W_q^I}=\|\phi\|_{L^q_I}+\|\nabla\nabla\phi\|_{L^q_I}+\left\|\frac\partial{\partial
t}\phi\right\|_{L^q_I}.
\end{align*}
The space $\mathcal W_q^I$ is a Sobolev-type space. It will play an
important part in our proof of Theorem~\ref{thm_parab_ex}.

Before proceeding, we need to introduce a family of atlases on $M$.
If $x\in M$ and $r>0$, suppose $B(x,r)$ is the open ball in $M$
centered at $x$ of radius $r$. Given $s>0$, consider an atlas
$\big(U_k^s,\bar\alpha_{s,k}\big)_{k=1}^{N(s)}$ on $M$ such that the
following requirements hold:
\begin{enumerate}
\item
The map $\bar\alpha_{s,k}$ is a diffeomorphism from $U_k^s$ to
$\mathbb R^n$ if $U_k^s$ lies in $M^\circ$ and a diffeomorphism from
$U_k^s$ to $\mathbb R^n_+$ if $U_k^s$ intersects $\partial M$.
\item
For every $k=1,\ldots,N(s)$, the domain $U_k^s$ equals
$B(x^{s,k},s)$ with $x^{s,k}\in M$. The map $\bar\alpha_{s,k}$ takes
$x^{s,k}$ to the origin.
\item
Given $x\in M$, there exists $k$ such that $x\in B(x^{s,k},\frac
s2)$ and the distance from $x$ to $M\setminus B(x^{s,k},\frac s2)$
is greater than $\tilde ss$. This requirement must hold for some
number $\tilde s>0$ independent of $s$.
\item
If $1\le k_1<\cdots<k_{N_0}\le N(s)$, then $U_{k_1}^s\cap\cdots\cap
U_{k_{N_0}}^s=\emptyset$. This must hold for some natural number
$N_0$ independent of~$s$.
\item
The inequality
\begin{align*}
\frac1{c'}|d\bar\alpha_{s,k}(\xi)|_{\mathbb R^n}\le|\xi|\le
c'|d\bar\alpha_{s,k}(\xi)|_{\mathbb R^n}
\end{align*}
is satisfied for all $\xi$ tangent to $M$ at a point $x$ whenever
$x\in B(x^{s,k},\frac{3s}4)$. Here, $|\cdot|_{\mathbb R_n}$ is the
standard Euclidean norm in $\mathbb R^n$, and $c'>0$ is a constant
independent of $s$, $k$, and $x$.
\end{enumerate}
It is clear that $\big(U_k^s,\bar\alpha_{s,k}\big)_{k=1}^{N(s)}$
exists as long as $s$ is sufficiently small;
cf.~\cite[page~295]{OLVSNU68}. In what follows, we fix $s_0>0$ such
that we can construct
$\big(U_k^s,\bar\alpha_{s,k}\big)_{k=1}^{N(s)}$ when $s\in(0,s_0]$.
It will be convenient for us to extend the diffeomorphisms
$\bar\alpha_{s,k}$ to the sets $U_k^s\times[0,I]$. More precisely,
we introduce the mappings $\alpha_{s,k}$ on $U_k^s\times[0,I]$ by
letting $\alpha_{s,k}(x,t)=(\bar\alpha_{s,k}(x),t)$.

Let us define a few cut-off functions. For each $s\in(0,s_0]$ and
$k=1,\ldots,N(s)$, choose a smooth $\bar\kappa_{s,k}:M\to[0,1]$
identically equal to~1 on $B(x^{s,k},\frac s2)$ and to~0 on
$M\setminus B(x^{s,k},\frac{3s}4)$. We may assume the norm of the
gradient of $\bar\kappa_{s,k}$ is bounded by $\frac{c''}s$ on $M$,
while the norm of the Hessian of $\bar\kappa_{s,k}$ is bounded by
$\frac{c''}{s^2}$ for some constant $c''>0$ independent of $s$ and
$k$. It will be convenient for us to define
\begin{align*}
\bar\rho_{s,k}=\frac{\bar\kappa_{s,k}}{\sum_{k=1}^{N(s)}(\bar\kappa_{s,k})^2}\,.
\end{align*} Finally, we introduce the function $\kappa_{s,k}$ on
$M\times[0,I]$ by setting $\kappa_{s,k}(x,t)=\bar\kappa_{s,k}(x)$.

Assume that, for every $s\in(0,s_0]$ and $k=1,\ldots,N(s)$, there is
a local trivialization $\bar\beta_{s,k}:\pi^{-1}(U_k^s)\to
U_k^s\times\mathbb R^d$ of the bundle $E$. This does not lead to any
loss of generality. Along with $\bar\beta_{s,k}$, let us introduce a
mapping $\beta_{s,k}:\pi^{-1}(U_k^s)\to\mathbb R^d$. By definition,
the image of $\theta\in\pi^{-1}(U_k^s)$ under $\beta_{s,k}$ is the
projection of $\bar\beta_{s,k}(\theta)$ onto $\mathbb R^d$. We
assume that the standard Euclidean norm of $\beta_{s,k}(\theta)$
equals $|\theta|$ for all $\theta\in\pi^{-1}(U_k^s)$ and, if $U_k^s$
intersects $\partial M$, the equality
\begin{align}\label{pr_descr1}
\beta_{s,k}(\pi^{-1}(x)\cap
W)&=\{(e_1,\ldots,e_{d'},0,\ldots,0)\in\mathbb
R^d\,|\,e_1,\ldots,e_{d'}\in\mathbb R\},\qquad x\in
U_k^s\cap\partial M,
\end{align}
holds for some $d'$ between~0 and~$d$. Again, these assumptions do
not lead to any loss of generality.

Let $(U_{k_l}^{s_0})_{l=1}^{N_1}$ be the collection of all those
$U_k^{s_0}$ that intersect $\partial M$. We use the notation $V_l$
for $U_{k_l}^{s_0}\cap\partial M$. It will be convenient for us to
write $\hat\alpha_l$ and $\hat\kappa_l$ for the restrictions of
$\bar\alpha_{s_0,k_l}$ and $\bar\kappa_{s_0,k_l}$ to $V_l$. One may
view $\hat\alpha_l$ as a diffeomorphism from $V_l$ to $\mathbb
R^{n-1}$. We can extend $\hat\alpha_l$ and $\hat\kappa_l$ to the
maps $\check\alpha_l$ and $\check\kappa_l$ on $V_l\times[0,I]$ by
setting $\check\alpha_l(x,t)=(\hat\alpha_l(x),t)$ and
$\check\kappa_l(x,t)=\hat\kappa_l(x)$. Also, let $\check\beta_l$ be
the restriction of $\beta_{s_0,k_l}$ to $\pi^{-1}(V_l)$.

Denote $\delta=1-\frac1q\,$, where $q$ is the integer fixed above.
We consider the fractional-order Sobolev-type space
$W^{\delta,\frac\delta2}_q(\mathbb R^{n-1}\times(0,I))$ of
real-valued functions. Its precise definition may be found
in~\cite[Chapter~II]{OLVSNU68}. Given a mapping $\phi:\partial
M\times[0,I]\to W^\bot$ such that $\pi(\phi(x,t))=x$, let us write
$\phi_l^m$ for the function taking $(y,t)\in\mathbb
R^{n-1}\times[0,I]$ to the $m$th component of the vector
$\big(\check\beta_l\circ\check\kappa_l\phi\circ\check\alpha_l^{-1}\big)(y,t)\in\mathbb
R^d$. Here, $l$ is an integer between $1$ and $N_1$, while $m$ is an
integer between $1$ and $d$. We introduce the space
$W^{\delta,\frac\delta2}_q(\partial M\times[0,I]\to W^\bot)$. It
consists of the mappings $\phi:\partial M\times[0,I]\to W^\bot$ such
that $\pi(\phi(x,t))=x$ for all $x\in\partial M$ and the function
$\phi_l^m$ lies in $W^{\delta,\frac\delta2}_q(\mathbb
R^{n-1}\times(0,I))$ for all $l=1,\ldots,N_1$ and $m=1,\ldots,d$.
The norm of $\phi$ in $W^{\delta,\frac\delta2}_q(\partial
M\times[0,I]\to W^\bot)$ is defined as the sum of the norms of
$\phi_l^m$ in $W_q^{\delta,\frac\delta2}(\mathbb
R^{n-1}\times(0,I))$. We use the notation
$\|\phi\|_{W^{\delta,\frac\delta2}_{q,I}}$ for it. Let us also
introduce the space $\mathcal W^{\delta,\frac\delta2}_{q,I}$ of
those $\phi\in W^{\delta,\frac\delta2}_q(\partial M\times[0,I]\to
W^\bot)$ that satisfy $\phi(x,0)=o(x)$ for all $x\in\partial M$.
Clearly, $\mathcal W^{\delta,\frac\delta2}_{q,I}$ inherits the norm
$\|\cdot\|_{W^{\delta,\frac\delta2}_{q,I}}$ from
$W^{\delta,\frac\delta2}_q(\partial M\times[0,I]\to W^\bot)$. The
nature of~$\mathcal W^{\delta,\frac\delta2}_{q,I}$ is explained in
part by the discussion on page~312 of~\cite{OLVSNU68}. Roughly
speaking, this space consists of the normal derivatives of the
mappings from $\mathcal W_q^I$.

\subsection{Linear parabolic equations}
\label{subsec_lin_par}

In order to prove Theorem~\ref{thm_parab_ex}, we need to establish a
few facts about second-order linear parabolic equations for sections
of $E$. We will heavily use material from~\cite{OLVSNU68}. Let us
lay down our setup.

Consider a smooth mapping
\begin{align*}
K:M\times[0,I]\to TM\otimes TM.
\end{align*}
Assume that $K(x,t)$ is a symmetric tensor over $x$ for all $x$ and
$t$. One may view $K$ as a section of the bundle $TM\otimes TM$
depending on $t\in[0,I]$. Consider one more mapping
\begin{align*}
G:M\times[0,I]&\to E.
\end{align*}
We demand that $G\in L^q_I$. Our interest is in the equation
\begin{align}\label{lin_par_eq}
v_t(x,t)-K^{ij}(x,t)\nabla_i\nabla_j v(x,t)=G(x,t),\qquad x\in
M^\circ,~t\in(0,I).
\end{align}
The unknown $v$ is a section of $E$ dependent on $t\in[0,I)$. The
subscript $t$ designates the differentiation in~$t\in(0,I)$.

Consider yet another mapping
\begin{align*}
p:\partial M\times[0,I]&\to W^\bot.
\end{align*}
We suppose it lies in $\mathcal W^{\delta,\frac\delta2}_{q,I}$. Let
us supplement~\eqref{lin_par_eq} with the boundary condition
\begin{align}\label{lin_par_BC}
\Pr_{W^\bot}\left(K^{ij}(x,t)\nu_i(x)\nabla_jv(x,t)\right)&=p(x,t),\qquad
x\in\partial M,~t\in(0,I).
\end{align}
Also, we assume there is a constant $c_2>0$ such that
\begin{align}\label{lin_par_cond}
K^{ij}(x,t)\xi_i\xi_j\ge c_2|\xi|^2
\end{align}
for every $x\in M$, $t\in[0,I]$, and $\xi\in T^*M$ projecting on
$x$. It is now time to state the main result of this subsection. It
gives us a solution to
problem~\eqref{lin_par_eq}--\eqref{lin_par_BC} in the space
$\mathcal W_q^I$ as well as an important estimate. The proof will be
largely based on the arguments in~\cite[Chapter~IV]{OLVSNU68}; see
also~\cite{VS65} and~\cite[Chapter~VII]{OLVSNU68}.

\begin{lemma}\label{lem_lin_par}
The boundary value problem~\eqref{lin_par_eq}--\eqref{lin_par_BC},
subject to condition~\eqref{lin_par_cond}, has a unique solution $v$
in the space $\mathcal W_q^I$. Furthermore, there exists $a>0$ such
that $v$ satisfies the estimate
\begin{align}\label{Schauder_est}
\|v\|_{\mathcal W_q^I}\le
a\Big(\|G\|_{L^q_I}+\|p\|_{W^{\delta,\frac\delta2}_{q,I}}\Big).
\end{align}
\end{lemma}

\begin{remark}
The number $a>0$ can be chosen independent of $I\in(0,1)$. Just how
large it has to be is determined by, among other things, the mapping
$K$.
\end{remark}

\begin{proof}[Proof of Lemma~\ref{lem_lin_par}]
We will produce the solution $v$ and establish~\eqref{Schauder_est}
assuming the number $I\in(0,1)$ is less than some number $I_0>0$ to
be specified later. In the end of the proof, we will remove this
assumption. Meanwhile, let $\mathfrak H_q^I$ stand for the direct
sum $L^q_I\oplus\mathcal W^{\delta,\frac\delta2}_{q,I}$. Define the
operator $\mathcal A:\mathcal W_q^I\to\mathfrak H_q^I$ by setting
$\mathcal Aw=(\mathcal A_1w,\mathcal A_2w\big)$ with
\begin{align*}
(\mathcal A_1w)(x,t)&=w_t(x,t)-K^{ij}(x,t)\nabla_i\nabla_j
w(x,t),\qquad x\in M^\circ,~t\in(0,I),
\\ (\mathcal A_2w)(x,t)&=\Pr_{W^\bot}\left(K^{ij}(x,t)\nu_i(x)\nabla_jw(x,t)\right),\qquad
x\in\partial M,~t\in(0,I).
\end{align*}
We need to show that $\mathcal A$ has a bounded inverse $\mathcal
A^{-1}$. The assertions of the lemma will follow immediately. The
role of the constant $a>0$ in inequality~\eqref{Schauder_est} will
be played by the norm of $\mathcal A^{-1}$.

To demonstrate that $\mathcal A$ has a bounded inverse, we blend the
arguments from the proofs of Theorems~5.3 and~9.1
in~\cite[Chapter~IV]{OLVSNU68}. Note that the geometric nature of
problem~\eqref{lin_par_eq}--\eqref{lin_par_BC} forces us to modify
those arguments rather substantially. Our first step is to
construct, assuming $I$ is less than $I_0$, a bounded operator
$\mathcal B:\mathfrak H_q^I\to\mathcal W_q^I$ such that the norms of
the operators $\mathcal A\mathcal B-\Id_{\mathfrak H_q^I}$ and
$\mathcal B\mathcal A-\Id_{\mathcal W_q^I}$ are less than~1. Here,
$\Id_{\mathfrak H_q^I}$ and $\Id_{\mathcal W_q^I}$ are the identity
maps in the corresponding spaces. Once $\mathcal B$ is at hand, we
will utilize it to produce a left inverse and a right inverse
for~$\mathcal A$. The existence of $\mathcal A^{-1}$ will be a
direct consequence.

Suppose $(J_1,J_2)\in\mathfrak H_q^I$. In order to specify how the
operator $\mathcal B$ acts on $(J_1,J_2)$, let us fix $s\in(0,s_0]$
and an atlas $(U_k^s,\bar\alpha_{s,k})_{k=1}^{N(s)}$ on $M$ as
described in Section~\ref{subsec_Holder_sp}. Along with
$(U_k^s,\bar\alpha_{s,k})_{k=1}^{N(s)}$, we have the collection
$(\beta_{s,k})_{k=1}^{N(s)}$. Each $\beta_{s,k}$ is a mapping from
$\pi^{-1}(U_k^s)$ to $\mathbb R^d$. Choose a domain $U_k^s$
intersecting $\partial M$. The diffeomorphism $\bar\alpha_{s,k}$
takes $U_k^s$ to $\mathbb R_+^n$. Let $\{y_1,\ldots,y_n\}$ be the
standard coordinates on $\mathbb R_+^n$. Given $y\in\mathbb R_+^n$
and $t\in[0,I]$, we write $\hat K^{ij}_{s,k}(y,t)$ and
$\hat\nu_i^{s,k}(y)$ for the components of the tensors
$K(\bar\alpha_{s,k}^{-1}(y),t)$ and $\nu(\bar\alpha_{s,k}^{-1}(y))$
with respect to $\bar\alpha_{s,k}$ and $\{y_1,\ldots,y_n\}$.

To describe the action of $\mathcal B$ on $(J_1,J_2)$, some
preparations are required. It will be convenient for us to denote
\begin{align*}
\hat J_{1,s,k}(y,t)&=\big(\beta_{s,k}\circ\kappa_{s,k}J_1\circ\alpha_{s,k}^{-1}\big)(y,t),\qquad y\in\mathbb R_+^n,~t\in(0,I),\\
\hat
J_{2,s,k}(y,t)&=\big(\beta_{s,k}\circ\kappa_{s,k}J_2\circ\alpha_{s,k}^{-1}\big)(y,t),\qquad
y\in\partial\mathbb R_+^n,~t\in(0,I).
\end{align*}
Consider the equation
\begin{align}\label{Eucl_eq}
z_t(y,t)-\hat K_{s,k}^{ij}(0,0)z_{y_iy_j}(y,t)=\hat J_{1,s,k}(y,t),
\qquad y\in\mathbb R^n_{+,0},~t\in(0,I).
\end{align}
Here, $z=(z^1,\ldots,z^d)$ is the vector of unknown functions with
$z^m:\mathbb R^n_+\times[0,I)\to\mathbb R$ for $m=1,\ldots,d$. The
subscript $y_i$ and $y_j$ mean component-wise differentiation in
$y_i$ and $y_j$. Let us impose the boundary conditions
\begin{align}\label{Eucl_BC}
\Pr_{\beta_{s,k}(\pi^{-1}(\bar\alpha_{s,k}^{-1}(y))\cap
W)}z(y,t)&=0, \notag \\
\Pr_{\beta_{s,k}(\pi^{-1}(\bar\alpha_{s,k}^{-1}(y))\cap
W^\bot)}\big(\hat
K_{s,k}^{ij}(0,0)\hat\nu_i^{s,k}(y)z_{y_j}(y,t)\big)&=\hat
J_{2,s,k}(y,t),\qquad y\in\partial\mathbb R^n_+,~t\in(0,I),
\end{align}
and the initial condition
\begin{align}\label{Eucl_IC}
z(y,0)=0,\qquad y\in\mathbb R^n_+.
\end{align}
Using formula~\eqref{pr_descr1}, we see that
problem~\eqref{Eucl_eq}--\eqref{Eucl_BC}--\eqref{Eucl_IC} is
equivalent to $d$ uncoupled problems, one for each $z^m$. The
Dirichlet boundary condition is imposed on $z^m$ when
$m=1,\ldots,d'$, and a Neumann-type condition is in place when
$m=d'+1,\ldots,d$. On the basis of Theorem~6.1
from~\cite[Chapter~IV]{OLVSNU68} (see also the argument on pages
343--345 of~\cite{OLVSNU68}), one easily concludes that
problem~\eqref{Eucl_eq}--\eqref{Eucl_BC}--\eqref{Eucl_IC} has a
unique solution whose components all lie in the space
$W^{2,1}_q(\mathbb R^n_{+,0}\times(0,I))$. We denote this solution
by $\hat{\mathcal B}^{s,k}_{J_1,J_2}$. Given $m=1,\ldots,d$, let
$\hat J_{1,s,k}^m$, $\hat J_{2,s,k}^m$, and $\hat{\mathcal
B}^{s,k,m}_{J_1,J_2}$ be the $m$th components of $\hat J_{1,s,k}$,
$\hat J_{2,s,k}$, and $\hat{\mathcal B}^{s,k}_{J_1,J_2}$. The
estimate
\begin{align}\label{est_Eucl}
\big\|\hat{\mathcal B}^{s,k,m}_{J_1,J_2}\big\|_{W^{2,1}_q(\mathbb
R^n_{+,0}\times(0,I))}\le a_1\Bigg(\bigg(\int_{\mathbb
R^n_+\times[0,I]}\big|\hat
J_{1,s,k}^m\big|^q\,dy\,dt\bigg)^{\frac1q} +\big\|\hat
J_{2,s,k}^m\big\|_{W^{\delta,\frac\delta2}_q(\mathbb
R^{n-1}\times(0,I))}\Bigg)
\end{align}
must be satisfied for some $a_1>0$ ($dy$ and $dt$ are the Lebesgue
measures on $\mathbb R_+^n$ and $[0,I]$). We introduce the map
$\mathcal B^{s,k}_{J_1,J_2}:U_k^s\times[0,I]\to E$ by setting
\begin{align*}
\mathcal
B^{s,k}_{J_1,J_2}(x,t)=\bar\beta_{s,k}^{-1}\big(x,\big(\hat{\mathcal
B}^{s,k}_{J_1,J_2}\circ\alpha_{s,k}\big)(x,t)\big), \qquad x\in
U_k^s,~t\in[0,I].
\end{align*}
It will be convenient for us to have $\mathcal B^{s,k}_{J_1,J_2}$
defined on all of $M\times[0,I]$. Therefore, we let $\mathcal
B^{s,k}_{J_1,J_2}(x,t)=o(x)$ for $x\in M\setminus U_k^s$ and
$t\in[0,I]$. The map $\mathcal B^{s,k}_{J_1,J_2}$ will be a
substantial ingredient in the image of $(J_1,J_2)$ under~$\mathcal
B$.

So far, we've been assuming $U_k^s$ intersected $\partial M$.
Suppose now $U_k^s$ is a domain contained in $M^\circ$. The
diffeomorphism $\bar\alpha_{s,k}$ then acts from $U_k^s$ to $\mathbb
R^n$. We can write down an equation analogous to~\eqref{Eucl_eq} in
$\mathbb R^n$ and an initial condition analogous to~\eqref{Eucl_IC}
in $\mathbb R^n$. The resulting problem will have a unique solution
with components in $W^{2,1}_q(\mathbb R^n\times(0,I))$. Again, we
denote this solution by $\hat{\mathcal B}^{s,k}_{J_1,J_2}$ and
introduce the map $\mathcal B^{s,k}_{J_1,J_2}$. Notice that each
component of $\hat{\mathcal B}^{s,k}_{J_1,J_2}$ will satisfy an
estimate similar to~\eqref{est_Eucl}. We are now ready to specify
how $\mathcal B$ acts on $(J_1,J_2)$.

Consider the mapping $\mathcal B_{J_1,J_2}^s:M\times[0,I]\to E$
given by the formula
\begin{align*}
\mathcal
B_{J_1,J_2}^s(x,t)=\sum_{k=1}^{N(s)}\bar\rho_{s,k}(x)\mathcal
B^{s,k}_{J_1,J_2}(x,t), \qquad x\in M,~t\in[0,I].
\end{align*}
It is clear that $\mathcal B_{J_1,J_2}^s\in\mathcal W_q^I$. Let
$\mathcal B$ be the operator taking $(J_1,J_2)\in\mathfrak H_q^I$ to
$\mathcal B_{J_1,J_2}^s\in\mathcal W_q^I$. Our next step is to
establish a few inequalities for $\mathcal B$.

In the beginning of the proof, we assumed $I$ was less than some
$I_0>0$. Fix a number $\bar s\in(0,1]$ and suppose $I_0$ is $\bar
ss^2$. Thus, the formula
\begin{align}\label{def_bar_s}
I<\bar ss^2
\end{align}
must hold. Employing~\eqref{est_Eucl}, \eqref{def_bar_s}, and the
properties of the atlas $(U_k^s,\bar\alpha_{s,k})_{k=1}^{N(s)}$
listed in Section~\ref{subsec_Holder_sp}, we easily see that the
inequality
\begin{align}\label{norm_cal_B}
\|\mathcal Bh\|_{\mathcal W_q^I}&<a_2\|h\|_{\mathfrak H_q^I},\qquad
h\in\mathfrak H_q^I,
\end{align}
is satisfied for some $a_2>0$ independent of the numbers $s$ and
$\bar s$; cf. Lemma~4.7 and Theorem~7.1
in~\cite[Chapter~IV]{OLVSNU68}. In particular, the operator
$\mathcal B$ is bounded.

If $s$ and $\bar s$ are chosen sufficiently small
and~\eqref{def_bar_s} holds, then
\begin{align*}
\|\mathcal A\mathcal Bh-h\|_{\mathfrak H_q^I}&<\|h\|_{\mathfrak H_q^I},\qquad h\in\mathfrak H_q^I,\\
\|\mathcal B\mathcal Aw-w\|_{\mathcal W_q^I}&<\|w\|_{\mathcal
W_q^I},\qquad w\in\mathcal W_q^I.
\end{align*}
In order to prove this, it is necessary to write down a series of
estimates based on~\eqref{norm_cal_B} and the H\"older inequality.
These estimates are very similar to the ones on pages~348-349
of~\cite{OLVSNU68}; see also~\cite{VS65}. We will not present them
here. At this point, the required inequalities for $\mathcal B$ are
at hand. We will now utilize $\mathcal B$ to produce the left
inverse and the right inverse of $\mathcal A$.

If $s$ and $\bar s$ are chosen small and~\eqref{def_bar_s} holds,
then the norms of the operators $\mathcal A\mathcal B-\Id_{\mathfrak
H_q^I}$ and $\mathcal B\mathcal A-\Id_{\mathcal W_q^I}$ are less
than~1. In this case, $\mathcal A\mathcal B=(\mathcal A\mathcal
B-\Id_{\mathfrak H_q^I})+\Id_{\mathfrak H_q^I}$ and $\mathcal
B\mathcal A=(\mathcal B\mathcal A-\Id_{\mathcal
W_q^I})+\Id_{\mathcal W_q^I}$ must have bounded inverses. Keeping
this in mind, we conclude that
\begin{align*}\mathcal A\left(\mathcal B(\mathcal A\mathcal
B)^{-1}\right)&=\Id_{\mathfrak H_q^I}, \\ \left((\mathcal B\mathcal
A)^{-1}\mathcal B\right)\mathcal A&=\Id_{\mathcal
W_q^I}.\end{align*} It is now easy to see that $\mathcal A$ must
have a bounded inverse $\mathcal A^{-1}$. Thus, the assertions of
the lemma hold true provided $I$ is less than $\bar ss^2$ for some
$s$ and $\bar s$. In order to complete the proof, we have to remove
the assumption on $I$. But this can be accomplished by repeating the
arguments from pages~349--350 of~\cite{OLVSNU68} (see
also~\cite[Chapter~IV, Section~8]{OLVSNU68}).
\end{proof}

\subsection{Proof of the existence theorem}
\label{subsec_par_pf}

Our preparations for the proof of Theorem~\ref{thm_parab_ex} are now
completed. We proceed in two steps. First, we will use a fixed-point
argument similar to the one found in~\cite{PW91} (see
also~\cite{MGGM87,PABT88}) to construct a solution $u$ of
problem~\eqref{parab_eq}--\eqref{parab_BC}--\eqref{parab_InitC}.
Then we will employ classical facts from~\cite{OLVSNU68} to show
that $u$ possesses the desired differentiability properties.

\begin{proof}[Proof of Theorem~\ref{thm_parab_ex}]
Let us assume $u_0(x)$ is equal to zero for every $x\in M$. This
does not lead to any loss of generality. Indeed, it is always
possible to reduce the general case to the case where $u_0(x)=o(x)$
for all $x\in M$ by introducing the new unknown $\check
u(x,t)=u(x,t)-u_0(x)$.

We will now construct the mapping $\mathcal C$ whose fixed point
will be a solution of
problem~\eqref{parab_eq}--\eqref{parab_BC}--\eqref{parab_InitC}. As
in Sections~\ref{subsec_Holder_sp} and~\ref{subsec_lin_par}, suppose
$I\in(0,1)$ is a real number and $q>n+2$ is an integer. The space
$\mathcal W_q^I$ will play an important role in our further
considerations. Denote $H_0(x)=H(o(x),0)$ for each $x\in M$. Suppose
that $w\in\mathcal W_q^I$. Let us introduce the mapping
$H_w:M\times[0,I]\to TM\otimes TM$ by the formula
\begin{align*}
H_w(x,t)&=H(w(x,t),t),\qquad x\in M,~t\in[0,I].
\end{align*}
The notation
\begin{align*}
F_{w,\nabla w}(x,t)&=F(w(x,t),\nabla w(x,t),t),\qquad x\in
M,~t\in[0,I],
\\ \Psi_w(x,t)&=\Psi(w(x,t),t),\qquad x\in\partial M,~t\in[0,I],
\end{align*}
will also be helpful. We consider the equation
\begin{align}\label{lin-quas_eq}
v_t(x,t)&-H_0^{ij}(x)\nabla_i\nabla_jv(x,t) \notag
\\ &=F_{w,\nabla
w}(x,t)+\big(H_w^{ij}(x,t)-H_0^{ij}(x)\big)\nabla_i\nabla_jw(x,t),\qquad
x\in M^\circ,~t\in(0,I),
\end{align}
for the unknown section $v$ depending on $t\in[0,I]$. We then
supplement this equation with the boundary condition
\begin{align}\label{lin-quas_bc}
\Pr_{W^\bot}&\big(H_0^{ij}(x)\nu_i(x)\nabla_jv(x,t)\big) \notag
\\ &=\Psi_w(x,t)-
\Pr_{W^\bot}\big(\big(H_w^{ij}(x,t)-H_0^{ij}(x)\big)\nu_i(x)\nabla_jw(x,t)\big),
\qquad x\in\partial M,~t\in(0,I).
\end{align}
Lemma~\ref{lem_lin_par} above demonstrates that
problem~\eqref{lin-quas_eq}--\eqref{lin-quas_bc} has a unique
solution in the space $\mathcal W_q^I$. We may, therefore, define a
mapping $\mathcal C:\mathcal W_q^I\to\mathcal W_q^I$ that takes $w$
to this solution. A series of estimates based
on~\eqref{Schauder_est} show the existence of a number $T\in(0,1)$
such that $\mathcal C$ has a fixed point when $I\le T$. These
estimates are very similar to the ones in the proofs of Lemmas~2.4
and~2.5 in~\cite{PW91}. We will not present them here. It is, thus,
possible to find $w\in\mathcal W_q^I$ satisfying the equality
$\mathcal C(w)=w$ provided $I\le T$. In particular, there exists
$u\in\mathcal W_q^T$ with $\mathcal C(u)=u$. Lemma~3.3
of~\cite[Chapter~II]{OLVSNU68} implies that $u$ and $\nabla u$ are
continuous on $M\times[0,T)$. Formulas~\eqref{parab_eq} and
\eqref{parab_BC} hold for $u$. We also have~\eqref{parab_InitC}
since we assumed $u_0(x)=o(x)$ for $x\in M$. It remains to show that
$u$ is smooth on $M\times(0,T)$. We will do so on the basis of
bootstrapping argument.

Fix an atlas $(U_k^s,\bar\alpha_{s,k})_{k=1}^{N(s)}$ on $M$ as
described in Section~\ref{subsec_Holder_sp}. Here, $s$ is an
arbitrary positive number less than $s_0$. Along with
$(U_k^s,\bar\alpha_{s,k})_{k=1}^{N(s)}$, we have the collection
$(\beta_{s,k})_{k=1}^{N(s)}$. Choose a domain $U_k^s$
intersecting~$\partial M$. The diffeomorphism $\bar\alpha_{s,k}$
takes $U_k^s$ to $\mathbb R_+^n$. We introduce the function $\tilde
u_{s,k}=\beta_{s,k}\circ u\circ\alpha_{s,k}^{-1}$. It acts from
$\mathbb R_+^n\times[0,T]$ to $\mathbb R^d$. We take $m=1,\ldots,d$
and write $\tilde u_{s,k}^m$ for the $m$th component of $\tilde
u_{s,k}$. Our next step is to study the differentiability of $\tilde
u_{s,k}^m$. This will help us obtain the desired conclusion about
the smoothness of $u$.

Let $\{y_1,\ldots,y_n\}$ be the standard coordinates in $\mathbb
R^n_+$. In what follows, the notation $D\tilde u_{s,k}$ stands for
the Jacobian matrix of $\tilde u_{s,k}$ with respect to
$\{y_1,\ldots,y_n\}$. It is not difficult to understand on the basis
of~\eqref{parab_eq} that $\tilde u_{s,k}^m$ satisfies
\begin{align}\label{aux_eq}
(\tilde u_{s,k}^m)_t(y,t)-\tilde H_{s,k,u}^{ij}(y,t)(\tilde
u_{s,k}^m)_{y_iy_j}(y,t)=\tilde F^{s,k,m}_{u,\nabla u}(y,t),\qquad
y\in\mathbb R^n_{+,0},~t\in(0,T).
\end{align}
In this formula,
\begin{align*}
\tilde H_{s,k,u}^{ij}(y,t)=\tilde H_{s,k}^{ij}(\tilde
u_{s,k}(y,t),y,t),\qquad y\in \mathbb R_{+,0}^n,~t\in(0,T),
\end{align*} with $\tilde H_{s,k}^{ij}$ being a function from $\mathbb R^d\times\mathbb R^n_+\times[0,T]$
to $\mathbb R$ for any $i,j=1,\ldots,n$. We point out that $\tilde
H_{s,k}^{ij}$ is smooth in all three variables. Also, in
formula~\eqref{aux_eq},
\begin{align*}
\tilde F_{u,\nabla u}^{s,k,m}(y,t)=\tilde F_{s,k}^m(\tilde
u_{s,k}(y,t), D\tilde u_{s,k}(y,t),y,t),\qquad y\in\mathbb
R_{+,0}^n,~t\in(0,T),
\end{align*}
with $\tilde F_{s,k}^m$ taking $\mathbb
R^d\times\Mat_{d,n}\times\mathbb R^n_+\times[0,T]$ to $\mathbb R$.
The notation $\Mat_{d,n}$ refers to the space of $d\times n$
matrices. We naturally identify $\Mat_{d,n}$ with $\mathbb R^{dn}$
and equip it with the standard Euclidean metric. Then the
function~$\tilde F_{s,k}^m$ is smooth in all four variables.
According to~\eqref{parab_BC} and~\eqref{pr_descr1}, if $m$ is
between~1 and $d'$, the equality
\begin{align}\label{aux_b1}
\tilde u_{s,k}^m(y,t)=0, \qquad y\in\partial\mathbb
R_+^n,~t\in(0,T),\end{align} holds true. For the other values of
$m$, we have a slightly more complicated identity. Finally, it is
easy to see from~\eqref{parab_InitC} that
\begin{align}\label{aux_i1}
\tilde u_{s,k}^m(y,0)=0, \qquad y\in\mathbb R_+^n.\end{align} We
proceed to establishing differentiability properties of $\tilde
u_{s,k}^m$.

The fact that $u$ lies in~$\mathcal W_q^T$ and Lemma~3.3
from~\cite[Chapter~II]{OLVSNU68} tell us that the coefficients
$\tilde H_{s,k,u}^{ij}$ and the term $\tilde F^{s,k,m}_{u,\nabla u}$
in equation~\eqref{aux_eq} lie in the H\"older-type space
$H^{\lambda,\frac\lambda2}(\Theta_{s,k}\times(0,T))$ for
$\Theta_{s,k}=\bar\alpha_{s,k}(B(x^{s,k},\frac s2))$ and some
$\lambda\in(0,1)$.
Taking~\eqref{aux_eq}--\eqref{aux_b1}--\eqref{aux_i1} into account
and employing the material in~\cite[Chapter~III,
Section~12]{OLVSNU68}, we can conclude that $\tilde u_{s,k}^m$ must
belong to $H^{2+\lambda,1+\frac\lambda2}(\Theta_{s,k}\times(0,T))$
if $m=1,\ldots,d'$. Analogous reasoning works when
$m=d'+1,\ldots,d$. In this case, $\tilde u_{s,k}^m$ is a solution
of~\eqref{aux_eq} under a Neumann-type boundary condition and the
initial condition~\eqref{aux_i1}. We easily see that $\tilde
u_{s,k}^m\in
H^{2+\lambda,1+\frac\lambda2}(\Theta_{s,k}\times(0,T))$.

Let us examine equation~\eqref{aux_eq} again. It is now evident that
the coefficients $\tilde H_{s,k,u}^{ij}$ and the term $\tilde
F^{s,k,m}_{u,\nabla u}$ in it must belong to
$H^{1+\lambda,\frac12+\frac\lambda2}(\Theta_{s,k}\times(0,T))$. As
above, we conclude $\tilde u_{s,k}^m\in
H^{3+\lambda,\frac32+\frac\lambda2}(\Theta_{s,k}\times(0,T))$. Let
us iterate this argument. It becomes clear that $\tilde u_{s,k}^m$
must be smooth on $\Theta_{s,k}\times(0,T)$.

So far, it's been assumed that the domain $U_k^s$ intersected
$\partial M$. Suppose now $U_k^s$ is contained in $M^\circ$. We can
still introduce the function $\tilde u_{s,k}=\beta_{s,k}\circ
u\circ\alpha_{s,k}^{-1}$. It now acts from $\mathbb R^n\times[0,T]$
to $\mathbb R^d$. Repeating the above reasoning with minor
modifications, we can demonstrate that, given $m=1,\ldots,d$, the
component $\tilde u_{s,k}^m$ must be smooth on
$\bar\alpha_{s,k}(B(x^{s,k},\frac s2))\times(0,T)$.

Let us summarize. Our goal was to establish the differentiability
properties of $u$ on $M\times(0,T)$. The presented arguments suggest
that, whichever $k=1,\ldots,N(s)$ we choose, $u$ must be smooth on
$B(x^{s,k},\frac s2)\times(0,T)$. This immediately implies the
desired properties of $u$ on $M\times(0,T)$.
\end{proof}

\section{The Ricci flow}\label{sec_RF}

Like in Section~\ref{sec_parab}, we consider a smooth manifold $M$
with boundary. We now assume that $M$ is $n$-dimensional with
$n\ge3$, compact, connected, and oriented. The notations $M^\circ$
and $\partial M$ will be used for the interior and the boundary of
$M$. Our goal is to investigate the Ricci flow on $M$. More
specifically, we will propose a new boundary condition and establish
two short-time existence results for this flow. The proofs will be
based on the method commonly known as DeTurck's trick. The reader
should consult, e.g.,~\cite{BCDK04,BCPLLN06,PT06} for a detailed
explanation of this method in the context of closed manifolds. A
relevant historic discussion may be found in~\cite{NCLG10}. Our
proofs will rely heavily on Theorem~\ref{thm_parab_ex}.

We focus on the equation
\begin{align}\label{RicciFlow_eq}
\frac{\partial}{\partial t}g(x,t)=-2\Ric^g(x,t),\qquad x\in
M^\circ,~t\in(0,T),
\end{align}
for a Riemannian metric $g$ on $M$ depending on the parameter
$t\in[0,T)$ with $T>0$. The notation $\Ric^g$ in the right-hand side
refers to the Ricci curvature of $g$. We fix a smooth Riemannian
metric $\hat g$ on $M$ and supplement~\eqref{RicciFlow_eq} with the
initial condition
\begin{align}\label{RicciFlow_IC}
g(x,0)=\hat g(x),\qquad x\in M.
\end{align}
So far, we do not concern ourselves with the behavior of $g$ near
$\partial M$. The reader will recognize that~\eqref{RicciFlow_eq} is
the Ricci flow equation on $M$. The introduction to the present
paper contains references to several books that discuss it in great
detail.

We call a mapping $g:M\times[0,T)\to T^*M\otimes T^*M$ a
\emph{decent solution} of
problem~\eqref{RicciFlow_eq}--\eqref{RicciFlow_IC} on $M\times[0,T)$
if the following requirements are met:
\begin{enumerate}
\item
For every $x\in M$ and $t\in[0,T)$, the tensor $g(x,t)$ is symmetric
and positive-definite. In other words, $g$ is a Riemannian metric on
$M$ depending on $t\in[0,T)$.
\item
The mapping $g$ is continuous on $M\times[0,T)$ and smooth on
$M\times(0,T)$.
\item
The Ricci flow equation~\eqref{RicciFlow_eq} and the initial
condition~\eqref{RicciFlow_IC} hold for $g$.
\end{enumerate}
Throughout Section~\ref{sec_RF}, we write $\tilde\nabla$ and
$\hat\nabla$ for the Levi-Civita connections of the metrics $g$ and
$\hat g$. In a similar fashion, $\tilde\upsilon$ and $\hat\upsilon$
will stand for the outward unit normal vector fields on $\partial M$
with respect to $g$ and $\hat g$. We point out that $\tilde\nabla$
and $\tilde\upsilon$ depend on the parameter $t\in[0,T)$, while
$\hat\nabla$ and $\hat\upsilon$ do not. The connection $\hat\nabla$
gives rise to connections in tensor bundles over $M$. We preserve
the notation $\hat\nabla$ for them.

As in Section~\ref{sec_parab}, let us implicitly assume that a
coordinate system $\{x_1,\ldots,x_n\}$ is chosen in a neighborhood
of every point $x\in M$. Suppose $T$ is a $(k,l)$-tensor field on
$M$ near $x$. By analogy with the notation of
Section~\ref{sec_parab}, we write $T_{j_1\ldots j_l}^{i_1\ldots
i_k}$ for the components of $T$ in the coordinates
$\{x_1,\ldots,x_n\}$, while $\tilde\nabla_iT$ and $\hat\nabla_i T$
stand for $\tilde\nabla_{\frac\partial{\partial x_i}}T$ and
$\hat\nabla_{\frac\partial{\partial x_i}}T$. The expression
$\hat\nabla_i\hat\nabla_jT$ means the second covariant derivative
$\hat\nabla\hat\nabla T$ applied to $\frac\partial{\partial x_i}$
and $\frac\partial{\partial x_j}$. If $x$ lies in $\partial M$, we
assume that $\{x_1,\ldots,x_{n-1}\}$ is a local coordinate system on
$\partial M$, the $n$th coordinate of any point in $\partial M$ near
$x$ is equal to~0, and $\hat\upsilon$ is a scalar multiple of
$\frac\partial{\partial x_n}$ near $x$. Given a $(k,l)$-tensor $Z$
on $\partial M$ at $x\in\partial M$, we write
$Z_{\beta_1\ldots\beta_l}^{\alpha_1\ldots\alpha_k}$ for the
components of $Z$ with respect to $\{x_1,\ldots,x_{n-1}\}$. As
before, the Einstein summation convention is in effect. The Latin
indices $i$, $j$, $k$, and $l$ will vary from 1 to $n$, whereas the
Greek indices $\alpha$, $\beta$, $\gamma$, and $\sigma$ will vary
from 1 to $n-1$.

In accordance with the notation introduced above, $g_{ij}$ and $\hat
g_{ij}$ are the components of the Riemannian metrics $g$ and $\hat
g$. We will also deal with the inverses of these metrics. Their
components will be denoted by~$g^{ij}$ and $\hat g^{ij}$.

\subsection{Formulation of the existence
results}\label{subsec_RF_results}

Our further considerations involve the second fundamental form field
$\II:\partial M\times[0,T)\to T^*\partial M\otimes T^*\partial M$ of
the boundary with respect to $g$. By definition,
\begin{align*}
\II_{\alpha\beta}(x,t)=g_{\alpha\gamma}(x,t)\big(\tilde\nabla_\beta\tilde\upsilon\big)^\gamma(x,t),\qquad
x\in\partial M,~t\in[0,T).
\end{align*}
Let us introduce the quantity
\begin{align*}
\mathcal
H(x,t)=\frac1{n-1}\,g^{\alpha\beta}(x,t)\II_{\alpha\beta}(x,t),\qquad
x\in\partial M,~t\in[0,T).
\end{align*}
It is called the mean curvature of $\partial M$. One may also
consider the second fundamental form field of $\partial M$ with
respect to $\hat g$. We will denote it by $\hat\II$. Finally, one
may introduce the mean curvature of $\partial M$ with respect
to~$\hat g$. We will write $\hat{\mathcal H}$ for it.

Let us state the first result of this section. It assumes that
$\hat{\mathcal H}(x)$ is independent of $x\in\partial M$. If this is
the case, we can solve
problem~\eqref{RicciFlow_eq}--\eqref{RicciFlow_IC} for a short time
maintaining control over $\mathcal H(x,t)$.

\begin{theorem}\label{thm_RF_mean_curv}
Suppose the mean curvature $\hat{\mathcal H}(x)$ is equal to a
constant $\mathcal H_0\in\mathbb R$ for all $x\in\partial M$. Let
$\mu$ be a smooth real-valued function on $[0,\infty)$ with
$\mu(0)=1$. Then there exist $T>0$ and a decent solution $g$ of
problem~\eqref{RicciFlow_eq}--\eqref{RicciFlow_IC} on $M\times[0,T)$
such that the mean curvature $\mathcal H(x,t)$ satisfies the
boundary condition
\begin{align*}
\mathcal H(x,t)=\mu(t)\mathcal H_0
\end{align*}
for all $x\in\partial M$ and $t\in[0,T)$.
\end{theorem}

\begin{remark}
We emphasize that the smooth function $\mu$ appearing in the theorem
can be arbitrary as long as $\mu(0)=1$. Essentially, different
choices of this function correspond to different evolutions of $\hat
g$ under the Ricci flow. It is also reasonable to think of $\mu$ as
a normalization factor. The number $T$ whose existence the theorem
asserts may depend on~$\mu$. The explicit form of this dependence,
however, is quite difficult to track down. We refer
to~\cite{MBXCAP10} for a discussion relevant to the geometric
meaning of $\mu$.
\end{remark}

The second result of this section touches upon the question of the
behavior of the Ricci flow on manifolds with convex boundary. Again,
it establishes the existence of a solution. Note that there are
several ways to define what it means for $\partial M$ to be convex
with respect to a Riemannian metric on $M$. Different viewpoints and
the relations between them are surveyed in~\cite{MS01}. Perhaps, the
most common way is to deem $\partial M$ convex with respect to a
Riemannian metric on $M$ if and only if the second fundamental form
of $\partial M$ with respect to this metric is nonnegative-definite
on $\partial M$. Having said that, we can formulate the result.

\begin{proposition}\label{thm_convex}
There exist $T>0$, a map $\psi$ from $M\times[0,T)$ to $M$, and a
decent solution $g$ of
problem~\eqref{RicciFlow_eq}--\eqref{RicciFlow_IC} on $M\times[0,T)$
such that the following statements hold:
\begin{enumerate}
\item\label{prop_1}
It is the case that $\psi$ is continuous on $M\times[0,T)$ and
smooth on $M\times(0,T)$.
\item\label{prop_2}
The map $\psi(\cdot,t)$ is a diffeomorphism from the manifold $M$ to
itself for every $t\in[0,T)$.
\item
The form field $\II(\cdot,t)$ coincides with the pullback of
$\hat{\II}$ by the restriction of $\psi(\cdot,t)$ to $\partial M$
whenever~$t\in[0,T)$.
\end{enumerate}
As a consequence, if $\hat\II(x)$ is nonnegative-definite for all
$x\in\partial M$ (that is, $\partial M$ is convex with respect to
$\hat g$), then $\II(x,t)$ is nonnegative-definite for all
$x\in\partial M$ and $t\in[0,T)$ (that is, $\partial M$ remains
convex with respect to $g$).
\end{proposition}

One more remark is in order at this point. After stating it, we will
proceed to proving Theorem~\ref{thm_RF_mean_curv} and
Proposition~\ref{thm_convex}. The first step will be to make some
preparations. We will do so in Section~\ref{subsec_Ricci-DeTurck}.

\begin{remark}
It may be possible to improve the regularity of $g$ on the set
$M\times[0,T)$ in Theorem~\ref{thm_RF_mean_curv} and
Proposition~\ref{thm_convex} by imposing additional restrictions on
the behavior of $\hat g$ near $\partial M$.
Remark~\ref{rem_par_comp} suggests the nature of the assumptions
that have to be made. We do not address this issue in the present
paper.
\end{remark}

\subsection{The Ricci-DeTurck flow and bundles over the boundary}\label{subsec_Ricci-DeTurck}

Let us introduce the Ricci-DeTurck flow on the manifold $M$. We will
employ it to prove Theorem~\ref{thm_RF_mean_curv}. Supplementing it
with boundary conditions will be the key step in our reasoning. The
existence of solutions to the flow will follow from
Theorem~\ref{thm_parab_ex}. We will then use a similar strategy to
prove Proposition~\ref{thm_convex}.

Given a Riemannian metric $\bar g$ on $M$ depending on $t\in[0,T)$,
define the mappings $P^{\bar g}:M\times[0,T)\to T^*M$ and $Q^{\bar
g}:M\times[0,T)\to T^*M\otimes T^*M$ by the formulas
\begin{align*}
P_i^{\bar g}(x,t)&=\bar g_{ij}(x,t)\bar g^{kl}(x,t)\big(\bar
\Gamma_{kl}^j(x,t)-\hat\Gamma_{kl}^j(x)\big),
\\
Q_{ij}^{\bar g}(x,t)&=\big(\bar\nabla_iP^{\bar
g}\big)_j(x,t)+\big(\bar\nabla_jP^{\bar g}\big)_i(x,t),\qquad x\in
M,~t\in[0,T).
\end{align*}
Here, $\bar g^{kl}$ are the components of the inverse of $\bar g$,
while $\bar \Gamma_{kl}^j$ and $\hat\Gamma_{kl}^j$ are the
Christoffel symbols corresponding to the Levi-Civita connections
$\bar\nabla$ and $\hat\nabla$ of $\bar g$ and $\hat g$. We have
written $\bar\nabla_iP^{\bar g}$ and $\bar\nabla_jP^{\bar g}$ for
the covariant derivatives of $P^{\bar g}$ in the directions
$\frac\partial{\partial x_i}$ and $\frac\partial{\partial x_j}$.
These derivatives are taken with respect to the connection in $T^*M$
induced by $\bar\nabla$. Consider the equation
\begin{align}\label{RDT_flow}
\frac{\partial}{\partial t}\bar g(x,t)=-2\Ric^{\bar g}(x,t)+Q^{\bar
g}(x,t),\qquad x\in M^\circ,~t\in(0,T).
\end{align}
In what follows, it will be convenient for us to use the notation
$\mathcal E$ for the bundle of symmetric (0,2)-tensors on $M$. One
could think about $\bar g$ as a section of $\mathcal E$ depending on
$t\in[0,T)$. By Lemma~7.48 in~\cite{BCPLLN06}, we can rewrite
equation~\eqref{RDT_flow} as
\begin{align}\label{RDT_fl_coord}
\frac{\partial}{\partial t}\bar g(x,t)=\bar
g^{ij}(x,t)\hat\nabla_i\hat\nabla_j \bar g(x,t)+R\big(\bar
g(x,t),\hat\nabla\bar g(x,t)\big), \qquad x\in M^\circ,~t\in(0,T).
\end{align}
The map $R$ here is defined on the set of all the pairs
$(\eta,\theta)$ such that $\eta$ is a symmetric positive-definite
(0,2)-tensor and $\theta$ is a (0,3)-tensor at the same point. The
values of $R$ lie in $T^*M\otimes T^*M$. Equation~\eqref{RDT_flow}
is the Ricci-DeTurck flow equation on $M$. Rewriting it in the
form~\eqref{RDT_fl_coord} will later enable us to apply
Theorem~\ref{thm_parab_ex} to it.

Let us introduce a few more pieces of notation. Suppose
$\pi_{\mathcal E}$ is the projection in the bundle $\mathcal E$. We
assume that $\mathcal E$ is equipped with the fiber metric given by
$\hat g$. It will be convenient for us to write $\mathcal
E_{\partial M}$ for the set of all $\eta\in\mathcal E$ such that
$\pi_{\mathcal E}(\eta)\in\partial M$. This set has the structure of
a vector bundle over $\partial M$ induced by the structure of
$\mathcal E$. It also inherits the fiber metric from $\mathcal E$.
Let $\mathcal F$ be the subbundle of $\mathcal E_{\partial M}$
consisting of all $\eta\in\mathcal E_{\partial M}$ such that
$\eta_{\alpha\beta}=0$ for $\alpha,\beta=1,\ldots,n-1$ and, in
addition, $\eta_{nn}=0$. One could view every $\eta\in\mathcal
E_{\partial M}$ as a bilinear form on $T_{\pi_{\mathcal E}(\eta)}M$.
With this interpretation adopted, the subbundle $\mathcal F$
consists of those $\eta\in\mathcal E_{\partial M}$ that satisfy
\begin{align*}
\eta(X,Y)=\eta(\hat\upsilon,\hat\upsilon)=0 \end{align*} for all $X$
and $Y$ tangent to $\partial M$ at $\pi_{\mathcal E}(\eta)$. Let
$\mathcal F^\bot$ be the orthogonal complement of $\mathcal F$ in
$\mathcal E_{\partial M}$.

\subsection{Proofs of the existence results}

We are now ready to prove Theorem~\ref{thm_RF_mean_curv}. In order
to do so, we will supplement the Ricci-DeTurck flow~\eqref{RDT_flow}
with boundary conditions and an initial condition.
Theorem~\ref{thm_parab_ex} will then imply the existence of a
solution. By modifying this solution, it is possible to obtain a
Riemannian metric on $M$ that
satisfies~\eqref{RicciFlow_eq}--\eqref{RicciFlow_IC} and exhibits
the desired boundary behavior.

\begin{proof}[Proof of Theorem~\ref{thm_RF_mean_curv}]
Recall that we are given a smooth real-valued function $\mu$ on
$[0,\infty)$ with $\mu(0)=1$. Let us impose boundary conditions on
the solutions of~\eqref{RDT_flow} by demanding that
\begin{align}\label{RDT_BC}
\Pr_{\mathcal F}\bar g(x,t)&=o(x), \notag \\
\bar\II_{\alpha\beta}(x,t)&=\frac12\,\mu(t)\big(\bar
g_{\alpha\gamma}(x,t)\hat
g^{\gamma\sigma}(x)\hat\II_{\sigma\beta}(x)+\bar
g_{\beta\gamma}(x,t)\hat
g^{\gamma\sigma}(x)\hat\II_{\sigma\alpha}(x)\big), \notag \\
P_n^{\bar g}(x,t)&=0,\qquad x\in\partial M,~t\in(0,T).
\end{align}
Here, $o$ is the zero section in $\mathcal E$, and $\bar\II$ is the
second fundamental form field of $\partial M$ with respect to $\bar
g$. A computation demonstrates that the boundary
conditions~\eqref{RDT_BC} are equivalent to the formulas
\begin{align}\label{RDT_boun_coord}
\Pr_{\mathcal F}\bar g(x,t)&=o(x), \notag \\
\Pr_{\mathcal F^\bot}\big(\bar g^{nn}(x,t)(\hat
g_{nn}(x))^\frac12\hat\nabla_n\bar g(x,t)\big)&=\zeta(\bar
g(x,t)),\qquad x\in\partial M,~t\in(0,T),
\end{align}
where $\zeta$ is a map from the set $\{\eta\in\mathcal E_{\partial
M}\,|\,\eta~\text{is positive-definite}\}$ to the bundle $\mathcal
F^\bot$. The components of the tensor $\zeta(\bar g(x,t))$ appear as
\begin{align*}
\zeta_{\alpha\beta}(\bar g(x,t))=&-\mu(t)\big(\hat g_{nn}(x)\bar
g^{nn}(x,t)\big)^{\frac12}\big(\bar g_{\alpha\gamma}(x,t)\hat
g^{\gamma\sigma}(x)\hat\II_{\sigma\beta}(x)+\bar
g_{\beta\gamma}(x,t)\hat
g^{\gamma\sigma}(x)\hat\II_{\sigma\alpha}(x)\big) \\
&+\hat g_{nn}(x)\bar g^{nn}(x,t)\big(\bar g_{\alpha\gamma}(x,t)\hat
g^{\gamma\sigma}(x)\hat\II_{\sigma\beta}(x)+\bar
g_{\beta\gamma}(x,t)\hat g^{\gamma\sigma}(x)\hat\II_{\sigma\alpha}(x)\big), \\
\zeta_{nn}(\bar g(x,t))= &-2\bar g_{nn}(x,t)\big(\mu(t)(\hat
g_{nn}(x)\bar g^{nn}(x,t))^{\frac12}\hat
g^{\alpha\beta}(x)\hat\II_{\alpha\beta}(x)-\bar
g^{\alpha\beta}(x,t)\hat\II_{\alpha\beta}(x)\big),
\end{align*}
and $\zeta_{\alpha n}(\bar g(x,t))=0$ whenever $x\in\partial M$ and
$t\in[0,T)$. As usual, the Greek indices vary from $1$ to $n-1$.

We supplement~\eqref{RDT_flow} with the initial condition
\begin{align}\label{RDT_IC}
\bar g(x,0)=\hat g(x),\qquad x\in M.
\end{align}
Theorem~\ref{thm_parab_ex} and Remark~\ref{rem_restr_HFg} imply the
existence of a number $T>0$ and a mapping $\bar
g:M\times[0,T)\to\mathcal E$ such that the following statements
hold:
\begin{enumerate}
\item The tensor $\bar g(x,t)$ is
positive-definite for every $x\in M$ and $t\in[0,T)$. In other
words, $\bar g$ is a Riemannian metric on $M$ depending on
$t\in[0,T)$.
\item
The mappings $\bar g$ and $\hat\nabla\bar g$ are continuous on
$M\times[0,T)$. Furthermore, $\bar g$ is smooth on $M\times(0,T)$.
\item Equalities~\eqref{RDT_fl_coord}, \eqref{RDT_boun_coord},
and~\eqref{RDT_IC} hold true.
\end{enumerate}
Evidently, $\bar g$ must also solve~\eqref{RDT_flow} and satisfy the
boundary conditions~\eqref{RDT_BC} along with the initial
condition~\eqref{RDT_IC}. We will write $\bar{\mathcal H}$ for the
mean curvature of $\partial M$ with respect to $\bar g$. Our next
step is to modify $\bar g$ by means of the DeTurck diffeomorphisms.
Then $\bar g$ will become a decent solution of
problem~\eqref{RicciFlow_eq}--\eqref{RicciFlow_IC}. Once we have
that, the proof of the theorem will be easy to complete.

Consider a mapping $P_{\bar g}:M\times[0,T)\to TM$ defined by the
formula
\begin{align*}
P_{\bar g}^i(x,t)=\bar g^{ij}(x,t)P^{\bar g}_j(x,t),\qquad x\in
M,~t\in[0,T).
\end{align*}
It is clear from~\eqref{RDT_BC} that $P_{\bar g}^n(x,t)$ must equal
0 when $x\in\partial M$ and $t\in[0,T)$. In other words, $P_{\bar
g}(x,t)$ is tangent to $\partial M$ as long as $x\in\partial M$. Let
us look at the equation
\begin{align}\label{ODE_DT}
\frac\partial{\partial t}\psi(x,t)=-P_{\bar g}(\psi(x,t),t),\qquad
x\in M,~t\in(0,T),
\end{align}
for $\psi:M\times[0,T)\to M$. We supplement this equation with the
initial condition
\begin{align}\label{ODE_IC}
\psi(x,0)=x,\qquad x\in M.
\end{align}
The mapping $P_{\bar g}$ is continuous on $M\times[0,T)$ and smooth
on $M\times(0,T)$. Also, $P_{\bar g}(\cdot,0)$ is identically zero
on $M$. Using these properties along with the fact that $P_{\bar
g}(x,t)$ is tangent to $\partial M$ whenever $x\in\partial M$, we
can prove the existence of a unique $\psi:M\times[0,T)\to M$ such
that
\begin{enumerate}
\item
The map $\psi$ is continuous on $M\times[0,T)$ and smooth on
$M\times(0,T)$.
\item
Equalities~\eqref{ODE_DT} and \eqref{ODE_IC} hold true.
\item
The map $\psi(\cdot,t)$ is a diffeomorphism from the manifold $M$ to
itself for every $t\in[0,T)$.
\end{enumerate}
The reader may find relevant material in~\cite[Chapter~3,
Section~3.1]{BCDK04}. It is customary to call $\psi(\cdot,t)$ the
DeTurck diffeomorphisms.

Given $t\in[0,T)$, define the Riemannian metric $g(\cdot,t)$ on $M$
as the pullback of $\bar g(\cdot,t)$ by $\psi(\cdot,t)$. One can
then verify that the mapping $g:M\times[0,T)\to T^*M\otimes T^*M$ is
a decent solution of
problem~\eqref{RicciFlow_eq}--\eqref{RicciFlow_IC} on the set
$M\times[0,T)$; see, e.g.,~\cite[page~81]{BCDK04}. Moreover, for
each $t\in[0,T)$, the second fundamental form field $\II(\cdot,t)$
is equal to the pullback of $\bar\II(\cdot,t)$ by the restriction of
$\psi(\cdot,t)$ to $\partial M$. Keeping this fact in mind,
taking~\eqref{RDT_BC} into account, and remembering the hypotheses
of the theorem, we compute the mean curvature $\mathcal H(x,t)$ and
see that
\begin{align*}
\mathcal H(x,t)=\bar{\mathcal H}(\psi(x,t),t)=\mu(t)\hat{\mathcal
H}(\psi(x,t))=\mu(t)\mathcal H_0,\qquad x\in\partial M,~t\in[0,T).
\end{align*}
The desired conclusion follows immediately.
\end{proof}

It is time to prove Proposition~\ref{thm_convex}. Again, we have to
supplement the Ricci-DeTurck flow~\eqref{RDT_flow} with boundary
conditions. The next step will be to apply
Theorem~\ref{thm_parab_ex} and obtain a solution. We will then
modify this solution by means of the DeTurck diffeomorphisms.

\begin{proof}[Proof of Proposition~\ref{thm_convex}]
Let us add boundary conditions to~\eqref{RDT_flow} by demanding that
\begin{align}\label{RDT_B2}
\Pr_{\mathcal F}\bar g(x,t)&=o(x), \notag \\
\bar\II_{\alpha\beta}(x,t)&=\hat\II_{\alpha\beta}(x),~P_n^{\bar
g}(x,t)=0,\qquad x\in\partial M,~t\in(0,T).
\end{align}
As before, $\bar\II$ is the second fundamental form field of
$\partial M$ with respect to $\bar g$. A computation shows that
formulas~\eqref{RDT_B2} are equivalent to
\begin{align*}
\Pr_{\mathcal F}\bar g(x,t)&=o(x), \notag \\
\Pr_{\mathcal F^\bot}\big(\bar g^{nn}(x,t)(\hat
g_{nn}(x))^\frac12\hat\nabla_n\bar g(x,t)\big)&=\chi(\bar
g(x,t)),\qquad x\in\partial M,~t\in(0,T),
\end{align*}
with $\chi$ acting from $\{\eta\in\mathcal E_{\partial
M}\,|\,\eta~\text{is positive-definite}\}$ to $\mathcal F^\bot$. The
components of $\chi(\bar g(x,t))$ are
\begin{align*}
\chi_{\alpha\beta}(\bar g(x,t))=&-2\big(\hat g_{nn}(x)\bar
g^{nn}(x,t)\big)^{\frac12}\hat\II_{\alpha\beta}(x) \\
&+\hat g_{nn}(x)\bar g^{nn}(x,t)\big(\bar g_{\alpha\gamma}(x,t)\hat
g^{\gamma\sigma}(x)\hat\II_{\sigma\beta}(x)+\bar
g_{\beta\gamma}(x,t)\hat g^{\gamma\sigma}(x)\hat\II_{\sigma\alpha}(x)\big), \\
\chi_{nn}(\bar g(x,t))=&-2\big((\hat g_{nn}(x)\bar
g_{nn}(x,t))^{\frac12}-\bar g_{nn}(x,t)\big)\bar
g^{\alpha\beta}(x,t)\hat\II_{\alpha\beta}(x),
\end{align*}
and $\chi_{\alpha n}(\bar g(x,t))=0$ when $x\in\partial M$ and
$t\in[0,T)$. As in the proof of Theorem~\ref{thm_RF_mean_curv}, we
apply Theorem~\ref{thm_parab_ex} and Remark~\ref{rem_restr_HFg} to
obtain $T>0$ and $\bar g:M\times[0,T)\to\mathcal E$ such that
\begin{enumerate}
\item The tensor $\bar g(x,t)$ is
positive-definite for every $x\in M$ and $t\in[0,T)$.
\item
The mappings $\bar g$ and $\hat\nabla\bar g$ are continuous on
$M\times[0,T)$. Also, $\bar g$ is smooth on $M\times(0,T)$.
\item
Equalities~\eqref{RDT_flow} and~\eqref{RDT_B2} hold true, and $\bar
g(x,0)=\hat g(x)$ whenever $x\in M$.
\end{enumerate}
It remains to bring the DeTurck diffeomorphisms into the picture and
to modify $\bar g$ by means of these diffeomorphisms. The latter
action will yield a decent solution
of~\eqref{RicciFlow_eq}--\eqref{RicciFlow_IC}.

Following the same procedure as in the proof of
Theorem~\ref{thm_RF_mean_curv}, let us define the mapping $P_{\bar
g}:M\times[0,T)\to TM$ and construct the corresponding
$\psi:M\times[0,T)\to M$. We introduce the Riemannian metric
$g(\cdot,t)$ on~$M$ as the pullback of $\bar g(\cdot,t)$ by
$\psi(\cdot,t)$ for every $t\in[0,T)$. One verifies that
$g:M\times[0,T)\to T^*M\otimes T^*M$ is a decent solution
of~\eqref{RicciFlow_eq}--\eqref{RicciFlow_IC} on $M\times[0,T)$.
Moreover, $\II(\cdot,t)$ coincides with the pullback of
$\bar\II(\cdot,t)$ by the restriction of $\psi(\cdot,t)$ to
$\partial M$ for each $t\in[0,T)$. This fact, along
with~\eqref{RDT_B2}, implies the desired conclusions.
\end{proof}

\section*{Acknowledgements}

I would like to thank Xiaodong Cao, Leonard Gross, and Carlos Kenig
for the productive discussions. I also express my gratitude to the
anonymous referee for his/her insightful review of the paper and
useful suggestions.


\begin{thebibliography}{99}

\bibitem{PABT88}
\textit{P. Acquistapace, B. Terreni}, On quasilinear parabolic
systems, Math. Ann.~\textbf{282} (1988), 315--335.

\bibitem{HA86}
\textit{H. Amann}, Quasilinear parabolic systems under nonlinear
boundary conditions,  Arch. Rational Mech. Anal.~\textbf{92} (1986),
153--192.

\bibitem{HA90}
\textit{H. Amann}, Dynamic theory of quasilinear parabolic equations
--- II. Reaction-diffusion systems, Diff. Int. Equ.~\textbf{3} (1990),
13--75.

\bibitem{MBXCAP10}
\textit{M. Bailesteanu, X. Cao, A. Pulemotov}, Gradient estimates
for the heat equation under the Ricci flow, J. Funct.
Anal.~\textbf{258} (2010), 3517--3542.

\bibitem{SB02}
\textit{S. Brendle}, Curvature flows on surfaces with boundary,
Math. Ann.~\textbf{324} (2002), 491--519.

\bibitem{NCLG10}
\textit{N. Charalambous, L. Gross}, The Yang-Mills heat semigroup on
three-manifolds with boundary, arXiv:1004.1639v1 [math.AP] (2010).

\bibitem{XCTD06}
\textit{X.-Z. Chen, T. Dong}, Ricci deformation of a metric on a
Riemannian manifold with boundary, J.~Zhejiang Univ. Sci.
Ed.~\textbf{33} (2006), 496--499.

\bibitem{BCDK04}
\textit{B. Chow, D. Knopf}, The Ricci flow: An introduction, Amer.
Math. Soc., Providence,~RI, 2004.

\bibitem{BCPLLN06}
\textit{B. Chow, P. Lu, L. Ni}, Hamilton's Ricci flow, Amer. Math.
Soc., Providence,~RI, 2006.

\bibitem{JC07}
\textit{J.C. Cortissoz}, The Ricci flow on the two ball with a
rotationally symmetric metric, Russian Math. (Iz. VUZ)~\textbf{51}
(2007), no. 12, 30--51.

\bibitem{JC09}
\textit{J.C. Cortissoz}, Three-manifolds of positive curvature and
convex weakly umbilic boundary, Geom. Dedicata~\textbf{138} (2009),
83--98.

\bibitem{MGGM87}
\textit{M. Giaquinta, G. Modica}, Local existence for quasilinear
parabolic systems under nonlinear boundary conditions, Ann. Mat.
Pura Appl.~\textbf{149} (1987), 41--59.

\bibitem{MHTW06}
\textit{M. Headrick, T. Wiseman}, Ricci flow and black holes, Class.
Quantum Grav.~\textbf{23} (2006), 6683--6707.

\bibitem{GHTSCW07}
\textit{G. Holzegel, T. Schmelzer, C. Warnick}, Ricci flows
connecting Taub-Bolt and Taub-NUT metrics, Class. Quantum
Grav.~\textbf{24} (2007), 6201--6217.

\bibitem{OLVSNU68}
\textit{O.A. Lady\v{z}enskaja, V.A. Solonnikov, N.N. Ural'ceva},
Linear and quasilinear equations of parabolic type, Amer. Math.
Soc., Providence,~RI, 1968.

\bibitem{TL93}
\textit{T. Li}, The Ricci flow on surfaces with boundary, Ph.D.
Dissertation, University of California at San Diego, 1993.

\bibitem{JMGT07}
\textit{J. Morgan, G. Tian}, Ricci Flow and the Poincar\'{e}
Conjecture, Amer. Math. Soc., Providence,~RI; Clay Mathematics
Institute, Cambridge,~MA, 2007.

\bibitem{TOVSEW06}
\textit{T.A. Oliynyk, V. Suneeta, E. Woolgar}, A gradient flow for
worldsheet nonlinear sigma models, Nuclear Phys.~B~\textbf{739}
(2006), 441--458.

\bibitem{TOEW07}
\textit{T.A. Oliynyk, E. Woolgar}, Rotationally symmetric Ricci flow
on asymptotically flat manifolds, Comm. Anal. Geom.~\textbf{15}
(2007), 535--568.

\bibitem{AP08}
\textit{A. Pulemotov}, The Li-Yau-Hamilton estimate and the
Yang-Mills heat equation on manifolds with boundary, J.~Funct.
Anal.~\textbf{255} (2008), 2933--2965.

\bibitem{MS01}
\textit{M. S\'{a}nchez}, Geodesic connectedness of semi-Riemannian
manifolds, Nonlinear Anal.~\textbf{47} (2001), 3085--3102.

\bibitem{YS92}
\textit{Y. Shen}, New results on some dynamical and stationary
problems in geometry, Ph.D.~Dissertation, Stanford University, 1992.

\bibitem{YS96}
\textit{Y. Shen}, On Ricci deformation of a Riemannian metric on
manifold with boundary, Pacific J.~Math.~\textbf{173} (1996),
203--221.

\bibitem{WXS89}
\textit{W.-X. Shi}, Ricci deformation of the metric on complete
noncompact Riemannian manifolds, J.~Differential Geom.~\textbf{30}
(1989), 303--394.

\bibitem{MS02}
\textit{M. Simon}, Deformation of $C\sp 0$ Riemannian metrics in the
direction of their Ricci curvature, Comm. Anal. Geom.~\textbf{10}
(2002), 1033--1074.

\bibitem{MS05}
\textit{M. Simon}, Deforming Lipschitz metrics into smooth metrics
while keeping their curvature operator non-negative, Geometric
evolution equations, Amer. Math. Soc., Providence, RI, 2005,
167--179.

\bibitem{VS65}
\textit{V.A. Solonnikov}, On boundary value problems for linear
parabolic systems of differential equations of general form, Proc.
Steklov Inst. Math.~\textbf{83} (1965), 1--184.

\bibitem{BT89}
\textit{B. Terreni}, Nonhomogeneous initial-boundary value problems
for linear parabolic systems, Studia Math.~\textbf{92} (1989),
141--175.

\bibitem{PT06}
\textit{P. Topping}, Lectures on the Ricci flow, Cambridge
University Press, Cambridge, 2006.

\bibitem{PW91}
\textit{P. Weidemaier}, Local existence for parabolic problems with
fully nonlinear boundary condition; an $L_p$-approach,  Ann. Mat.
Pura Appl.~\textbf{160}  (1991), 207--222.

\end{thebibliography}
\end{document}